\documentclass{article}

\usepackage[ansinew]{inputenc}
\usepackage{amsmath}
\usepackage{amssymb}
\usepackage{amsthm}
\usepackage{amscd}
\usepackage{amsfonts}
\usepackage[english]{babel}
\usepackage{graphicx}
\newtheorem {defn}{Definition}

\newtheorem {thm}[defn]{Theorem}

\newtheorem {lemma}[defn]{Lemma}

\newcommand{\I}{\boldsymbol{\mathcal{I}}}

\newcommand{\SU}{\boldsymbol{\mathcal{SU}}}
\newcommand{\U}{\boldsymbol{\mathcal{U}}}
\newcommand{\R}{\boldsymbol{\mathcal{R}}}

\newcommand{\RSU}{\boldsymbol{r\mathcal{SU}}}
\newcommand{\RU}{\boldsymbol{r\mathcal{U}}}
\newcommand{\RR}{\boldsymbol{r\mathcal{R}}}

\newcommand{\TR}{\operatorname{tr}}

\newcommand{\bbmatrix}[1]{\left[ \begin{array}{ccccccccc} #1 \end{array} \right]}
\newcommand{\SPEC}{{\operatorname{spec}}}

\newcommand{\Rem}{\noindent\textbf{Remark: }}
\newcommand{\Rems}{\noindent\textbf{Remarks: 1.}}
\newcommand{\Proof}{\noindent\textbf{Proof: }}

\title{Solution of the Boussinesq equation using evolutionary vessels}
\author{Andrey Melnikov\\Drexel University}

\begin{document}
\maketitle
\abstract{In this work we present a solution of the Boussinesq equation. The derived formulas include solitons, Schwartz class solutions and solutions,
possessing singularities on a closed set $Z$ of $\mathbb R^2$ ($(x,t)$ domain), obtained from the zeros of the tau function. 
The idea for solving the Boussinesq equation is identical to 
the (unified) idea of solving the KdV and the evolutionary NLS equations: we use a theory of evolutionary vessels. But
a more powerful theory of non-symmetric evolutionary vessels is presented, inserting flexibility into the construction
and allowing to deal with complex-valued solutions.
A powerful scattering theory of Deift-Tomei-Trubowitz for a three dimensional operator,
which is used to solve the Boussinesq equation, fits into our setting only in a particular case.
On the other hand, we create a much wider class of solutions of the Boussinesq equation with singularities 
on a closed set $Z$.}
\tableofcontents

\section{Introduction}
The Boussinesq equation \cite{bib:Bous}\footnote{The full Boussinesq equation $q_{tt} = \dfrac{\partial^2}{\partial x^2} [q + q_{xx} - 4q^2 ]$ was shown by McKean
\cite{bib:McKeanBoussinesq} to be equivalent to \eqref{eq:Boussinesq}, significantly simplifying algebra computations.}
\begin{equation} \label{eq:Boussinesq}
q_{tt} = \dfrac{\partial^2}{\partial x^2} [3q_{xx} - 12 q^2 ]
\end{equation}
is a foundation of a shallow water theory. A fundamental solution of this equation was presented by P. {D}eift, C. {T}omei and E. {T}rubowitz in \cite{bib:DTTBoussinesq}. The basic
idea in this work, presented by V. E. Zakharov \cite{bib:ZachForBouss}, is a development of a scattering theory of a three dimensional operator
\[ \widetilde L = i \dfrac{d^3}{dx^3} + \dfrac{1}{i} (q \dfrac{d}{dx} + \dfrac{d}{dx} q) + p
\]
and using the fact that for $Q = i (3 \dfrac{d^2}{dx^2} - 4 q)$, the operators $\widetilde L, Q$ constitute a Lax pair:
\[ \dfrac{d}{dt} \widetilde L = [Q,\widetilde L] = Q \widetilde L - \widetilde L Q.
\]
It is assumed in \cite{bib:DTTBoussinesq} that $p(x), q(x)$ are in a Schwartz space, but the theory goes through for $q(x), p(x)$ with only a finite number of derivatives
and a finite order of decay. The scattering data in this case consists of a list of 6 functions, whose evolving with $t$ under an analogue of $Q$ is studied producing
a solution of the Boussinesq equation \eqref{eq:Boussinesq} with a given initial value $q(x,0)=q(x)$. Since we will use instead of $\widetilde L$ its multiplication on $-i$, we
define
\begin{equation} \label{eq:DefL}
L = -i \widetilde L = \dfrac{d^3}{dx^3}  - 2 q \dfrac{d}{dx} - (q_x + i p)
\end{equation}
and a very particular (inverse) scattering of this operator will be researched. More precisely, we will discuss solutions of
\footnote{\label{ftn:derivs}we will usually denote by $q'$ the partial derivative $\dfrac{\partial}{\partial x}q$, and by $\dot q$ the partial derivative $\dfrac{\partial}{\partial t}q$. Similarly,
for the higher derivatives: $q''$ stands for $\dfrac{\partial^2}{\partial x^2}q$, etc. }
\begin{equation} \label{eq:Lisk3}
L u = \dfrac{d^3}{dx^3}  - 2 q \dfrac{d}{dx} - (q' + i p) = k^3 u, \quad p'(x)=P(\int q,q,q',\ldots,q^{(n)}).
\end{equation}
In this setting $q(x)$ is an ``arbitrary'' function, \textit{but $p(x)$ is derived from $q(x)$ (see formula \eqref{eq:Dp})} using the formula
$p'(x) = P(\int q,q,q',\ldots,q^{(n)})$ for a polynomial $P$, which is not derived explicitly, since we do not use it.

In order to solve \eqref{eq:Boussinesq}, we present a similar to \cite{bib:DTTBoussinesq} scheme,
where we use a special case of the inverse scattering theory: we create inverse scattering of \eqref{eq:DefL},
where $p(x)$ is not arbitrary and then evolve $q(x)$ with $t$. Let us explain first the (inverse) scattering theory.
The scattering data is encoded in our setting in a matrix-valued function.
In fact, there is an almost ``one-to-one''
correspondence (Theorem \ref{thm:Uniqueness}) between the coefficients $q(x)$, defining $L$ \eqref{eq:DefL} and $3\times 3$ matrix-functions $S(\lambda)$ possessing
a realization \cite{bib:bgr}
\[ S(\lambda) = I - C_0 \mathbb X_0^{-1} (\lambda I - A)^{-1} B_0 \sigma_1, \quad \sigma_1 = \bbmatrix{0 & 0 & 1 \\ 0 & 1 & 0 \\ 1 & 0 & 0}, \quad
A \mathbb X_0 + \mathbb X_0 A_\zeta + B_0\sigma_1C_0 = 0.
\]
Here for an auxiliary Hilbert space $\mathcal H$ the linear operators act as follows: $C_0:\mathcal H\rightarrow\mathbb C^3$, $A_\zeta,\mathbb X_0, A: \mathcal H\rightarrow\mathcal H$, $B_0:\mathbb C^3\rightarrow\mathcal H$. Let us assume for the simplicity of the presentation, that all the operators are bounded. In order to construct
$q(x)$, uniquely defined from $S(\lambda)$, use the following 4 steps, fixing
$\sigma_2 = \bbmatrix{1 & 0 & 0 \\ 0 & 0 & 0 \\ 0 & 0 & 0}$, $\gamma = \bbmatrix{0 & 0 & 0 \\ 0 & 0 & 1 \\ 0 & -1 & 0}$:
\begin{enumerate}
	\item solve for $B(x)$ from $B'\sigma_1= -AB\sigma_2-B\gamma$ \eqref{eq:DB} with initial $B(x_0)=B_0$,
	\item solve for $C(x)$ from $\sigma_1C' = - \sigma_2 C A_\zeta + C \gamma$ \eqref{eq:DC} with initial $C(x_0)=C_0$,
	\item solve for $\mathbb X(x)$ from $\mathbb X'(x) = B(x)\sigma_2C(x)$ \eqref{eq:DX}, $\mathbb X(x_0)=\mathbb X_0$,
	\item define $\gamma_*(x) = \gamma + \sigma_2C(x)\mathbb X^{-1}(x)B(x)\sigma_1 - \sigma_1C(x)\mathbb X^{-1}(x)B(x)\sigma_2$ \eqref{eq:Linkage},
	for all points where $\mathbb X(x)$ is invertible.
\end{enumerate}
Then we prove in Theorem \ref{thm:Backlund} that the function
\[ S(\lambda,x) = I - C(x) \mathbb X^{-1}(x) (\lambda I - A)^{-1} B(x) \sigma_1
\]
is a B\" acklund transformation for the operator $L$ from the trivial $L_0 = \dfrac{d^3}{dx^3}$ to a more complicated one $L$ \eqref{eq:DefL}, in which $q(x) = - \dfrac{3}{2}  \dfrac{d^2}{dx^2} \ln \det\big(\mathbb X_0^{-1}\mathbb X(x)\big)$ and $p(x)$ is defined from $q(x)$ up to a constant. In the regular case, explained here
the coefficients $q(x), p(x)$ are analytic functions at all points, where $\mathbb X(x)$ is invertible.

Letting the operators $C,B,\mathbb X$ further evolve with respect to $t$, 
we will obtain a solution of the Boussinesq equation
\eqref{eq:Boussinesq}. In order to show it, we take three additional matrices \eqref{eq:3DimCanParamt}
\[ \widetilde{\sigma}_1 = \sigma_1, \quad \widetilde\sigma_2=\bbmatrix{0&-i&0\\i&0&0\\0&0&0},\quad \widetilde\gamma=\bbmatrix{0&0&0\\0&0&0\\0&0&i}
\]
and follow these steps (last step is the same as the previous fourth step):
\begin{enumerate}
	\item[$\dot 4$.] solve for $B(x,t)$ from $\dot B \widetilde\sigma_1 = -AB\widetilde\sigma_2-B\widetilde\gamma$ \eqref{eq:DBt} (see footnote \ref{ftn:derivs}) with initial $B(x,t_0)=B(x)$,
	\item[$\dot 5$.] solve for $C(x,t)$ from $\widetilde\sigma_1\dot C = - \widetilde\sigma_2 C A_\zeta + C \widetilde\gamma$ \eqref{eq:DCt} with initial $C(x,t_0)=C(x)$,
	\item[$\dot 6$.] solve for $\mathbb X(x,t)$ from $\dot{\mathbb X}(x) = B(x,t)\widetilde\sigma_2C(x,t)$ \eqref{eq:DX}, $\mathbb X(x,t_0)=\mathbb X(x)$,
	\item[7.=4.] define $\gamma_*(x,t) = \gamma + \sigma_2C(x,t)\mathbb X^{-1}(x,t)B(x,t)\sigma_1 - \sigma_1C(x,t)\mathbb X^{-1}(x,t)B(x,t)\sigma_2$ \eqref{eq:Linkage},
	for all points where $\mathbb X (x,t)$ is invertible.
\end{enumerate}
It turns out that $q(x,t) = - \dfrac{3}{2}  \dfrac{\partial^2}{\partial x^2} \ln \det\big(\mathbb X_0^{-1}\mathbb X(x,t)\big)$ \eqref{eq:qtau}
and satisfies the Boussinesq equation \eqref{eq:Boussinesq} (see Theorem \ref{thm:BoussnesqProof} for details).
Construction of the coefficient $q(x)$ from a realized function $S(\lambda)$ is called \textit{the standard construction
of a vessel} and is presented in Section \ref{sec:GenVessels}.

The simplicity and richness of this construction is best revealed in soliton formulas (Section \ref{sec:solitons}).
By choosing the inner space $\mathcal H=\mathbb C$ - the one dimensional Hilbert space, we create a classical
soliton $q(x,t) = -\dfrac{9\mu^2}{2\cosh^2(\dfrac{\sqrt{3}}{2} \mu(x+t\mu))}$ \eqref{eq:classSol} and another one
$q(x,t) = -\dfrac{18e^{2\sqrt{3}\mu(x+2t\mu)}\mu^2}{(e^{2\sqrt{3}x\mu} + e^{4\sqrt{3}t\mu^2})^2}$ \eqref{eq:solnonclass}.
Here $\mu$ is an arbitrary complex parameter.

This construction seems to be more concrete and suitable for the solution of the Boussinesq equation \eqref{eq:Boussinesq}, because it uses
just ``enough'' of the very powerful and complicated inverse scattering theory, developed in \cite{bib:DTTBoussinesq}. Notice that in this later work
the coefficients $q(x), p(x)$ are quite arbitrary. In our work, on the other hand, the coefficient $p(x)$ is uniquely determined (up to a constant) from $q(x)$. Still, the formulas enable to produce solutions, applying the standard construction of a vessel, explained above, to different $S(\lambda)$. 
Moreover, at the same section \ref{sec:GenVessels} we show that
much more general classes of solutions arise, as we impose as few as possible restrictions on $S(\lambda)$. 
The solutions, presented in \cite{bib:DTTBoussinesq} are either in the Schwartz class or exponentially decaying
and correspond in our setting to an analytic
$S(\lambda)$, possessing jumps along the real negative axis.

Finally, this work presents a very general setting for a construction of solutions of \eqref{eq:Boussinesq}. We find
necessary regularity assumptions \eqref{eq:BregBouss} on the operator $B(x,t)$ such that we can create a vessel,
and hence a solution of \eqref{eq:Boussinesq}. This is proved in Theorem \ref{thm:BoussSolGen}.
Following the remark after Theorem \ref{thm:BoussSolGen}, one can easily construction a solution, 
which fails to be five times $x$-differentiable (four times are necessary for the existence of \eqref{eq:Boussinesq}).

\section{Scattering theory of the operator $L$}
We start from the definition of the vessel parameters, which create an inverse scattering theory of $L$ \eqref{eq:DefL}.
\begin{defn} \label{def:3DimCanParam}
The vessel parameters are defined as follows
\[ \sigma_1 = \bbmatrix{0 & 0 & 1 \\ 0 & 1 & 0 \\ 1 & 0 & 0},
   \sigma_2 = \bbmatrix{1 & 0 & 0 \\ 0 & 0 & 0 \\ 0 & 0 & 0},
   \gamma   = \bbmatrix{0 & 0 & 0 \\ 0 & 0 & 1 \\ 0 & -1 & 0}.
\]
\end{defn}

\subsection{Non-symmetric regular vessel, realizing a scattering theory of $L$}
\begin{defn} 
A (regular, non-symmetric) vessel, associated to vessel parameters (see Definition \ref{def:3DimCanParam}) is a collection of operators, spaces and an interval $\mathrm I$
\begin{equation} \label{eq:DefV}
\mathfrak{V}_{reg} = (C(x), A_\zeta, \mathbb X(x), A, B(x); \sigma_1, 
\sigma_2, \gamma, \gamma_*(x);
\mathcal{H},\mathbb C^3;\mathrm I),
\end{equation}
where the bounded operators $C(x):\mathcal H\rightarrow\mathbb C^3$, $A_\zeta, \mathbb X(x), A :\mathcal H\rightarrow\mathcal H$, $B(x):\mathbb C^3\rightarrow\mathcal H$ 
and a $3\times 3$ matrix function $\gamma_*(x)$ satisfy the following vessel conditions:
\begin{eqnarray}
\label{eq:DB} \frac{\partial}{\partial x} B  &  & = - (A \, B \sigma_2 + B \gamma)\sigma_1^{-1}, \\
\label{eq:DC} \frac{\partial}{\partial x} C  &  & = \sigma_1^{-1} (\gamma C - \sigma_2 C A_\zeta ), \\
\label{eq:DX} \frac{\partial}{\partial x} \mathbb X &  & =  B \sigma_2 C, \\
\label{eq:Linkage}
&  \gamma_*  & =  \gamma + \sigma_2 C \mathbb X^{-1} B \sigma_1 -
 \sigma_1 C \mathbb X^{-1} B \sigma_2, \\
\label{eq:Lyapunov} & A \, \mathbb X + \mathbb X\,  A_\zeta  & = - B \sigma_1 C.
\end{eqnarray}
The operator $\mathbb X(x)$ is assumed to be invertible on the interval $\mathrm I$. If $A_\zeta=A^*$ and $C=B^*$ we call such a vessel \textbf{symmetric}.
\end{defn}
\noindent\textbf{Remarks: 1.} Notice that the operators $C(x), \mathbb X(x), B(x)$ are globally defined for all $x\in\mathbb R$: $C(x), B(x)$ are solutions of 
operator-valued differential equations with constant coefficients and $\mathbb X(x)$ is obtained from them by a simple integration. 
\textbf{2.} For the definition of the matrix-function $\gamma_*(x)$ we need the invertability of the operator $\mathbb X(x)$, so in general we may suppose that $\gamma_*(x)$ is defined
for all $x\in\mathbb R$, except for those points where $\mathbb X(x)$ is not invertible. For simplicity, we take an interval $\mathrm I$ and in a more general setting
we consider $\gamma_*(x)$ on $\mathbb R$, except fro the points where $\mathbb X(x)$ is not invertible.
\begin{thm}[permanency of the Lyapunov equation] \label{thm:LyapPerm}
Suppose that $B(x), C(x), \mathbb X(x)$ satisfy  \eqref{eq:DB}, \eqref{eq:DC} and \eqref{eq:DX} respectively,
then if the Lyapunov equation \eqref{eq:Lyapunov}
\[ A \mathbb X(x) + \mathbb X(x) A_\zeta  + B(x) \sigma_1 C(x) = 0\]
holds for a fixed $x_0\in\mathrm I$, then it holds for all $x\in\mathrm I$.
\end{thm}
\noindent\textbf{Proof:} By differentiating the left hand side of \eqref{eq:Lyapunov}, we will obtain that 
it is zero. \qed

By the definition, the \textit{transfer function} of the vessel $\mathfrak{V}_{reg}$ is defined as follows:
\begin{equation}\label{eq:DefS} 
S(\lambda,x) = I - C(x) \mathbb X^{-1}(x) (\lambda I - A)^{-1} B(x) \sigma_1.
\end{equation}
Notice that poles and singularities of $S$ with respect to $\lambda$ are determined by $A$ only.
We would like to show that the function $S(\lambda,x)$ realizes a B\" acklund transformation of the corresponding LDEs:
multiplication by the function $S(\lambda,x)$ maps \cite{bib:Vortices, bib:SchurIEOT, bib:SLVessels} a solution of the input Linear Differential Equation
(LDE) with the spectral parameter $\lambda$
\begin{equation}
\label{eq:InCC}
 \lambda \sigma_2 u(\lambda, x) -
\sigma_1 \frac{\partial}{\partial x}u(\lambda,x) +
\gamma u(\lambda,x) = 0.
\end{equation}
to a solution of the output LDE with the same spectral parameter
\begin{equation} \label{eq:OutCC}
\lambda \sigma_2 y(\lambda, x) - \sigma_1 \frac{\partial}{\partial x}y(\lambda,x) +
\gamma_*(x) y(\lambda,x) = 0.
\end{equation}
The function $\gamma_*$ is defined by the Linkage condition \eqref{eq:Linkage}. The fundamental solutions of \eqref{eq:InCC} and \eqref{eq:OutCC}, which are equal to $I$ 
(the identity matrix matrix at $x=0$) are denoted usually by $\Phi(\lambda,x)$ and $\Phi_*(\lambda,x)$. In other words, these are matrix functions
satisfying:
\begin{eqnarray}
\label{eq:DefPhi} \lambda \sigma_2 \Phi(\lambda,x) -
\sigma_1 \frac{\partial}{\partial x}\Phi(\lambda,x) +
\gamma \Phi(\lambda,x) = 0, \quad \Phi(\lambda,0) = I, \\
\label{eq:DefPhi*} \lambda \sigma_2 \Phi_*(\lambda,x) -
\sigma_1 \frac{\partial}{\partial x}\Phi_*(\lambda,x) +
\gamma_*(x) \Phi_*(\lambda,x) = 0, \quad \Phi_*(\lambda,0) = I.
\end{eqnarray} 
\begin{thm}[Vessel=B\" acklund transformation \cite{bib:Vortices, bib:SchurIEOT, bib:SLVessels}] \label{thm:Backlund}Suppose that $u(\lambda, x)$ satisfies \eqref{eq:InCC}, then
$y(\lambda,x) = S(\lambda,x)u(\lambda,x)$ satisfies \eqref{eq:OutCC}.
\end{thm}
\noindent\textbf{Remark:} Notice that it is enough to prove that the transfer function satisfies the following ODE:
\begin{equation}\label{eq:DS} 
\dfrac{\partial}{\partial x} S(\lambda,x) =
\sigma_1^{-1}(\sigma_2\lambda+\gamma_*(x)) S(\lambda,x) -
S(\lambda,x) \sigma_1^{-1}(\sigma_2\lambda+\gamma).
\end{equation}
\noindent\textbf{Proof:} First we calculate
\[ \begin{array}{llll}
\dfrac{d}{dx} \big( C(x)\mathbb X^{-1}(x)\big) & = \sigma_1^{-1} (\gamma C -\sigma_2 C A_\zeta) \mathbb X^{-1} -
C \mathbb X^{-1} B \sigma_2 C \mathbb X^{-1} = \text{ by \eqref{eq:Lyapunov}} \\
& = \sigma_1^{-1} \gamma C\mathbb X^{-1} + \sigma_1^{-1} \sigma_2 C \mathbb X^{-1} A +  \sigma_1^{-1} \sigma_2 C \mathbb X^{-1} B\sigma_1 C \mathbb X^{-1}-
C \mathbb X^{-1} B \sigma_2 C \mathbb X^{-1} \\
& = \sigma_1^{-1} \sigma_2 C \mathbb X^{-1} A + \sigma_1^{-1} (\gamma + \sigma_2 C \mathbb X^{-1} B \sigma_1 -
 \sigma_1 C \mathbb X^{-1} B \sigma_2) C\mathbb X^{-1} = \text{ by \eqref{eq:Linkage}} \\
& =  \sigma_1^{-1} \sigma_2 C \mathbb X^{-1} A + \sigma_1^{-1} \gamma_* C\mathbb X^{-1}.
\end{array} \]
So we obtain that
\begin{equation} \label{eq:DCX}
\dfrac{d}{dx} \big( C\mathbb X^{-1} \big) = \sigma_1^{-1} \sigma_2 C \mathbb X^{-1} A + \sigma_1^{-1} \gamma_* C\mathbb X^{-1}.
\end{equation}
Let us differentiate next the transfer function using \eqref{eq:DefS}:
\[ \begin{array}{llll}
\dfrac{d}{dx} S(\lambda,x) & = - \dfrac{d}{dx} \big( C\mathbb X^{-1}\big) (\lambda I - A)^{-1} B\sigma_1 -
C\mathbb X^{-1}(\lambda I - A)^{-1}\dfrac{d}{dx}B\sigma_1  = \text{by \eqref{eq:DB}, \eqref{eq:DCX}} \\
& = \big( \sigma_1^{-1} \sigma_2 C \mathbb X^{-1} A + \sigma_1^{-1} \gamma_* C\mathbb X^{-1}\big)(\lambda I - A)^{-1} B\sigma_1 -
C\mathbb X^{-1}(\lambda I - A)^{-1}(A \, B \sigma_2 + B \gamma)\\
& =  \sigma_1^{-1} \gamma_* ( S- I) - (S - I)\sigma_1^{-1} \gamma + \sigma_1^{-1} \sigma_2 C \mathbb X^{-1} A (\lambda I - A)^{-1} B\sigma_1-
C\mathbb X^{-1}(\lambda I - A)^{-1} A \, B \sigma_2 \\
& = \text{insert $A=A\pm\lambda I$ and expand} \\
& = \sigma_1^{-1} \gamma_* ( S- I) - (S - I)\sigma_1^{-1} \gamma + \sigma_1^{-1} \sigma_2 C \mathbb X^{-1} B\sigma_1-
 \lambda \sigma_1^{-1} \sigma_2 C \mathbb X^{-1} (\lambda I - A)^{-1} B\sigma_1 -\\
 & \quad \quad + C\mathbb X^{-1}B \sigma_2 -\lambda C\mathbb X^{-1}(\lambda I - A)^{-1} B \sigma_2  \\
 & = \text{using \eqref{eq:Linkage}, \eqref{eq:DefS}} \\
 & = \sigma_1^{-1} \gamma_* ( S- I) - (S - I)\sigma_1^{-1} \gamma + \sigma_1^{-1}(\gamma_*-\gamma)+
   \lambda \sigma_1^{-1} \sigma_2 S - S  \sigma_1^{-1} \lambda \sigma_2 \\
 & = \sigma_1^{-1}(\lambda\sigma_2 + \gamma_*) S - S\sigma_1^{-1}(\lambda\sigma_2 + \gamma). 
\end{array} \]
\qed

Expanding the transfer function $S(\lambda,x)$ into a Taylour series around $\lambda=\infty$, we obtain a notion of the moment:
\[ S(\lambda,x) = I - \sum_{n=0}^\infty \dfrac{H_n(x)\sigma_1}{\lambda^{n+1}},
\]
where by the definition \textit{the $n$-th moment} $H_n(x)$ of the function $S(\lambda,x)$ is
\begin{equation} \label{eq:Hn}
H_n(x) = C(x)\mathbb X^{-1}(x)A^nB(x).
\end{equation}
Using the zero moment, for example, we obtain that the so called ``linkage condition'' \eqref{eq:Linkage} is equivalent to
\[ \gamma_*(x) = \gamma + \sigma_2 H_0(x)\sigma_1 - \sigma_1 H_0(x)\sigma_2.
\]
There is also a recurrent relation between the moments $H_n(x)$, arising from \eqref{eq:DS}:
\begin{thm} The following recurrent relation between the moments of the vessel $\mathfrak{V}_{reg}$ holds
\begin{eqnarray}
\label{eq:DHn} 	\sigma_1^{-1}\sigma_2H_{n+1} - H_{n+1} \sigma_2\sigma_1^{-1} = (H_n)'_x - \sigma_1^{-1} \gamma_* H_n + H_n \gamma\sigma_1^{-1}.
\end{eqnarray}
\end{thm}
\noindent\textbf{Proof: } Follows from the differential equation \eqref{eq:DS} by plugging $S(\lambda,x) = I - \sum_{n=0}^\infty \dfrac{H_n(x)\sigma_1}{\lambda^{n+1}}$. \qed

Let us investigate more carefully the LDEs \eqref{eq:InCC} and \eqref{eq:OutCC}. Denote $u = \bbmatrix{u_1\\u_2\\u_3}$, then
\eqref{eq:InCC} becomes
\[ \bbmatrix{\lambda u_1 \\0 \\0} - \bbmatrix{u_3'\\u_2'\\u_1'}+
\bbmatrix{0 \\  u_3
\\ - u_2 } = \bbmatrix{0\\0\\0}.
\]
Solving this we obtain that 
\begin{equation} \label{eq:uentries}\left\{ \begin{array}{lll}
u_2 = -u_1', \\
u_3 = u_2' = -u_1'', \\
u_1''' =- \lambda u_1. \\
\end{array}\right.
\end{equation}
We can see that actually this equation is equivalent to a third-order differential equation with the spectral parameter $\lambda$:
\begin{equation} \label{eq:InCCequiv}
u_1''' = -\lambda u_1.
\end{equation}
In order to analyze \eqref{eq:OutCC}, we denote first moment $H_0(x)= [\pi_{ij}] = \bbmatrix{ \pi_{11} & \pi_{12} & \pi_{13} \\ \pi_{21} & \pi_{22} & \pi_{23} \\
\pi_{31} & \pi_{32} & \pi_{33} }$, and as a result, the linkage condition \eqref{eq:Linkage} becomes
\[ \gamma_* = \gamma + \bbmatrix{\pi_{13}-\pi_{31} & \pi_{12}& \pi_{11} \\
-\pi_{21} & 0 & 0 \\ -\pi_{11} & 0 & 0}.
\]
Denote next $y=\bbmatrix{y_1\\y_2\\y_3}$ and plugging the expression for $\gamma_*$ just derived into \eqref{eq:OutCC}, we will obtain that
\[ \bbmatrix{\lambda y_1 \\0 \\0} - \bbmatrix{y_3'\\y_2'\\y_1'}+
\bbmatrix{(\pi_{13}-\pi_{31}) y_1 + \pi_{12} y_2 + \pi_{11}  y_3 \\ -\pi_{21}  y_1 + y_3
\\ -\pi_{11}  y_1 - y_2 } = \bbmatrix{0\\0\\0}\]
or solving this:
\begin{equation} \label{eq:entries}
\left\{ \begin{array}{lll}
y_2 = -\pi_{11}y_1 - y_1', \\
y_3 = \pi_{21}y_1 + y_2' = \pi_{21} y_1 - y_1 \pi_{11}' - 
 \pi_{11} y_1' - y_1'', \\
y_1''' - 2 q y_1' - (q' + p) y_1 = - \lambda y_1,
\end{array}\right.
\end{equation}
where
\begin{equation} \label{eq:Defpq}
q(x) = \dfrac{\pi_{11}^2 + \pi_{12}+\pi_{21} - 2 \pi_{11}'}{2}, \quad
p(x) = -i(-\pi_{13} + \pi_{31} +  \pi_{11} (\pi_{12}-\pi_{21}) - \dfrac{\pi_{12}'-\pi_{21}'}{2}).
\end{equation}
In other words the equation \eqref{eq:OutCC} is equivalent to
\begin{equation} \label{eq:outCCequiv}
y_1''' - 2 q y_1' - (q' + i p) y_1 = - \lambda y_1.
\end{equation}
Now we are ready to justify the term ``scattering matrix'' attached to $S(\lambda,0)$. An independent set 
\[ y(\lambda,x)= \bbmatrix{y_1(\lambda,x)\\y_2(\lambda,x)\\y_3(\lambda,x)}
\]
of solutions of \eqref{eq:DefL}, can be derived from \eqref{eq:OutCC} in the following form
\[ y(\lambda,x) = S(\lambda,x) \Phi(\lambda,x) \bbmatrix{1\\0\\0} \left(= \Phi_*(\lambda,x) S(\lambda,0)\bbmatrix{1\\0\\0}\right).
\]
The fundamental matrix $\Phi(\lambda,x)$, solving \eqref{eq:InCC} with the initial condition $\Phi(\lambda,0)=I$ (-the identity matrix) is obtained from
\eqref{eq:uentries} and can be explicitly written as
\begin{eqnarray}
\label{eq:PhiForm}\Phi(\lambda,x) = \dfrac{1}{3\alpha^2} \bbmatrix{
R_1 & \dfrac{R_2}{k} & \dfrac{-R_3}{k^2} \\
k R_3 & R_1 & -\dfrac{R_2}{k} \\
-k^2 R_2 & -k R_3 & R_1
} \\
\nonumber R_1 = \alpha^2(E_1 + E_2 + E_3), \\
\nonumber R_2 = \alpha^2 E_1 + \alpha E_2 + E_3, \\
\nonumber R_3 = \alpha^2E_1 + E_2 +  \alpha E_3
\end{eqnarray}
where $E_1=e^{-kx}, E_2=e^{-\alpha k x}, E_3=e^{-\alpha^2k x}$ for $\alpha = e^{2\pi i/3}$ ($\alpha^3=1$).
This matrix is analytic in $\lambda$, because examining Taylor series of $E_i$ we will come to the conclusion that all the entries
of $\Phi$ depend on $k^3=\lambda$. The structure of $S(\lambda,x)$ is also known from \eqref{eq:DefS}, so we can study solutions of \eqref{eq:DefL} or equivalently
of \eqref{eq:OutCC}, creating in this manner the (inverse) scattering of $L$ \eqref{eq:DefL}.

\subsection{Structure of the moment $H_0(x)$}
Let us examine the recurrence relation \eqref{eq:DHn}. We will research for the simplicity of the
presentation the structure of the first moment $H_0(x)$, but almost the same structure will actually 
apply for all moments. Let us denote
\[ H_1(x) = \bbmatrix{ g_{11} & g_{12} & g_{13} \\ g_{21} & g_{22} & g_{23} \\
g_{31} & g_{32} & g_{33} }.
\]
Then the left hand side of \eqref{eq:DHn} becomes
\[ \sigma_1^{-1}\sigma_2H_1(x) - H_1(x)\sigma_2\sigma_1^{-1} = 
\bbmatrix{0&0&-g_{11}\\0&0&-g_{21}\\g_{11}&g_{12}&g_{13}-g_{31}}
\]
which must be equal to
\[ \begin{array}{ll}
(H_0)'_x  - \sigma_1^{-1} \gamma_* H_0 + H_0 \gamma\sigma_1^{-1} = \\
= \bbmatrix
{\pi_{11}^2 +\pi_{12}+\pi_{21} + \pi_{11}'& \pi_{11} \pi_{12} - \pi_{13} + \pi_{22} + \pi_{12}'& \pi_{11}\pi_{13} + \pi_{23} + \pi_{13}'\\
\pi_{11} \pi_{21} - \pi_{31} + \pi_{22} + \pi_{21}' & \pi_{12}\pi_{21} - \pi_{23} -\pi_{32} + \pi_{22}' & Res(23) \\
-\pi_{11} \pi_{13} - \pi_{12}\pi_{21} + \pi_{32} + \pi_{31}'& Res(32) & Res(33)} ,
\end{array} \]
where
\[ \begin{array}{ll}
Res(23) = \pi_{13} \pi_{21} - \pi_{33} + \pi_{23}', \\
Res(32) = -\pi_{12} (\pi_{13} - \pi_{31}+ \pi_{22}) - \pi_{11} \pi_{32} - \pi_{33} +  \pi_{32}', \\
Res(33) = -\pi_{12} \pi_{23} + \pi_{13}(-\pi_{13} + \pi_{31})  - \pi_{11}\pi_{33} +  \pi_{33}'.
\end{array} \]
Equating we obtain that 
\begin{equation} \label{eq:knowngij}
g_{11} = - \pi_{11}\pi_{13} - \pi_{23} - \pi_{13}',\quad
g_{21}=-Res(23), \quad g_{12}=Res(32), \quad g_{13}-g_{31} = Res(33)
\end{equation}
and the following system of equations:
\[
 \left\{ \begin{array}{llllll}
\pi_{11}^2 + \pi_{12} + \pi_{21} + \pi_{11}' = 0, \\
\pi_{11} \pi_{12} - \pi_{13} + \pi_{22} + \pi_{12}' = 0, \\
\pi_{11} \pi_{21} - \pi_{31} + \pi_{22} + \pi_{21}' = 0, \\
\pi_{12}\pi_{21} - (\pi_{23}+\pi_{32}) + \pi_{22}' = 0, \\
- \pi_{11}\pi_{13} - \pi_{23} - \pi_{13}' = -\pi_{11} \pi_{13} - \pi_{12}\pi_{21} + \pi_{32} + \pi_{31}' (= g_{11}).
\end{array} \right. \]
or rearranging
\begin{equation} \label{eq:PIsimple}
\left\{ \begin{array}{llllll}
\pi_{21} = - (\pi_{12} + \pi_{11}^2 + \pi_{11}' ), \\
\pi_{22} = -\pi_{11} \pi_{12} + \pi_{13} - \pi_{12}', \\
\pi_{22} = -\pi_{11} \pi_{21} + \pi_{31} - \pi_{21}', \\
\pi_{32} = \pi_{12}\pi_{21} - \pi_{23} + \pi_{22}', \\
\pi_{32} = \pi_{12}\pi_{21} - \pi_{23} - \pi_{31}' - \pi_{13}'.
\end{array} \right. \end{equation}
Plugging the fourth equation of \eqref{eq:PIsimple} into the last one, we will obtain that the later becomes $\pi_{31}' + \pi_{13}' = -\pi_{22}'$ or requiring a
\textit{normalization}
\begin{equation} \label{eq:Normal1} \pi_{31}(x_0) + \pi_{13}(x_0)= -\pi_{22}(x_0)
\end{equation}
we obtain
\begin{equation} \label{eq:pi13pi22} \pi_{31} + \pi_{13}= -\pi_{22}. \end{equation}
In a similar manner, one can solve for some other entires and we summarize these intermediate results in the next lemma. Additional relations between the
entries of $H_0(x)$ are obtained if we consider \eqref{eq:DHn} for $n=1$. For example, formulas similar to \eqref{eq:PIsimple} are as follows:
\[
 \left\{ \begin{array}{llllll}
\pi_{11}g_{11} + g_{12} + g_{21} + g_{11}' = 0, \\
\pi_{11} g_{12} - g_{13} + g_{22} + g_{12}' = 0, \\
g_{11} \pi_{21} - g_{31} + g_{22} + g_{21}' = 0, \\
g_{12}\pi_{21} - (g_{23}+g_{32}) + g_{22}' = 0, \\
- \pi_{11}g_{13} - g_{23} - g_{13}' = -\pi_{11} g_{31} + g_{11} (\pi_{31}- \pi_{13})- \pi_{12}g_{21} + g_{32} + g_{31}'.
\end{array} \right. \]
From these formulas follow the folllowing relations
\begin{equation} \label{eq:gsimple}
 \left\{ \begin{array}{llllll}
\pi_{11}g_{11} + g_{12} + g_{21} + g_{11}' = 0, \\
\pi_{11} g_{12} - \pi_{11}g_{21} - (g_{13}-g_{31}) + g_{12}'-g'_{21} = 0,
\end{array} \right. \end{equation}
which can be rewritten as relations on the entries of $H_0(x)$ in view of \eqref{eq:knowngij}. Other three relations, which are obtained are
\[ \left\{ \begin{array}{llllll}
\pi_{11} g_{12} - g_{13} + g_{22} + g_{12}' + g_{11} \pi_{21} - g_{31} + g_{22} + g_{21}' = 0, \\
g_{12}\pi_{21} - (g_{23}+g_{32}) + g_{22}' = 0, \\
- \pi_{11}g_{13} - g_{23} - g_{13}' = -\pi_{11} g_{31} + g_{11} (\pi_{31}- \pi_{13})- \pi_{12}g_{21} + g_{32} + g_{31}',
\end{array} \right. \]
which serve to find $g_{22}'$, $g_{23} + g_{32}$, $g'_{13} + g'_{31}$. It turns out that actually the relations \eqref{eq:knowngij} are not independent and a formula for
$\pi_{12}\pi_{11}'$ is derived from them. All these results are summarized in the next Lemma and we notice that the exact
formulas for $\pi_{23}', \pi_{33}'$ are omitted, because we are not interested in their form:
\begin{lemma} The following relations between the entries of the first moment $H_0(x)$ hold
\begin{eqnarray}
\label{eq:pi21} \pi_{21} & = - (\pi_{12} + \pi_{11}^2 + \pi_{11}' ), \\
\label{eq:pi22} \pi_{22} & = -\pi_{11} \pi_{12} + \pi_{13} - \pi_{12}', \\
\label{eq:pi31} \pi_{31} & = \pi_{13} - \pi_{11}^3-\pi_{11}(2\pi_{12}+3\pi_{11}')-2\pi_{12}'-\pi_{11}'', \\
\label{eq:pi32} \pi_{32} & = \pi_{12}\pi_{21} - \pi_{23} -\pi_{12}\pi_{11}'-\dfrac{3}{2}((\pi_{11}')^2-\pi_{11}\pi_{11}'') + \dfrac{1}{2}\pi_{11}''', \\
\label{eq:pi13'} \pi_{13}' &  = -\dfrac{3}{2} (\pi_{11}')^2 + \pi_{11}\pi_{12}' + \dfrac{3}{2}\pi_{11}\pi_{11}''+\pi_{12}''+\dfrac{1}{2}\pi_{11}''', \\
\label{eq:pi12relation} 6 \pi_{12}\pi_{11}' & = -(6\pi_{11}+15\pi_{11}')\pi_{11}' + 3\pi_{11}\pi_{11}''+\pi_{11}'''.
\end{eqnarray}
The relations \eqref{eq:gsimple} result in formulas for $\pi_{23}', \pi_{33}'$.
\end{lemma}
\Proof Notice that we have obtained 5 relations on the entries of $H_0(x)$ in \eqref{eq:PIsimple}, 2 relations in \eqref{eq:gsimple} and one more relation
from the fact that \eqref{eq:knowngij} together with \eqref{eq:gsimple} are overdetermined. 
It is possible to derive these formulas using any symbolic computation software.
\qed

It would be interesting to derive formulas for the entries of $H_0(x)$, relying on the structure of the matrices $\sigma_1,\sigma_2,\gamma$ only,
without referring to tedious direct computations. It can be probably done using the following notion, appearing first in \cite{bib:SLVessels}:
\begin{defn} Tau-function of the vessel $\mathfrak V$ is defined as follows:
\begin{equation} \label{eq:Deftau}\tau(x) = \det(\mathbb X^{-1}(x_0)\mathbb X(x)),
\end{equation}
where $x_0\in\mathrm I$ is an arbitrary point.
\end{defn}
The fact that this object is well-defined follows from an equivalent to \eqref{eq:DX} equation
\[ \mathbb X(x) = \mathbb X_0 + \int_0^x B(y)\sigma_2C(y)dy,
\]
so that $\mathbb X_0^{-1}\mathbb X(x) = I + T(x)$ for a trace-class operator $T$.
Using this notion  we obtain the following
\begin{thm} The coefficient $q(x)$ possesses the following formula:
\begin{eqnarray}
\label{eq:qtau} q(x) & = - \dfrac{3}{2} \dfrac{d^2}{dx^2} \ln \tau(x) =  - \dfrac{3}{2} \dfrac{d}{dx} \pi_{11}.
\end{eqnarray}
The derivative of the coefficient $p(x)$ is a differential polynomial in $\pi_{11}$:
\begin{equation} \label{eq:Dp}
p'(x) = P(\pi_{11},\pi'_{11},\ldots), \quad \text{ $P$ - a polynomial}.
\end{equation}
\end{thm}
\noindent\textbf{Proof:} Inserting \eqref{eq:pi21} into the formula \eqref{eq:Defpq} for $q(x)$ we find that
\[ q(x) = \dfrac{\pi_{11}^2 + \pi_{12} +\pi_{21} - 2 \pi_{11}'}{2} =
-\dfrac{3}{2} \pi_{11}'.
\]
Using a formula for the determinant of an operator, we obtain that (see \cite{bib:GKintro, bib:SchurIEOT} for details)
\[ \dfrac{\tau'(x)}{\tau(x)} = \TR(\mathbb X'(x)\mathbb X^{-1}(x)) = \TR (B(x)\sigma_2B^*(x)\mathbb X^{-1}(x))=
\TR(\sigma_2H_0(x))=\pi_{11}
\]
and the result follows for $q(x)$. As for $p(x)$, we differentiate the formula for $p(x)$, appearing in \eqref{eq:Defpq}:
\[ p'(x) = -i(-\pi'_{13} + \pi'_{31}) -i \dfrac{d}{dx}( \pi_{11} (\pi_{12}-\pi_{21}) - \dfrac{\pi_{12}'-\pi_{21}'}{2}).
\]
Then using \eqref{eq:pi31}, \eqref{eq:pi13'}, \eqref{eq:pi21} and \eqref{eq:pi12relation} we will obtain a differential polynomial in $\pi_{11}$.
\qed

One can also derive the following formulas, corresponding to the symmetric case $C=B^*$, $A_\zeta=A^*$, plugging the definition \eqref{eq:Deftau} and
using the fact that $H_0(x)$ is self-adjoint:
\begin{lemma} For the symmetric case, $C=B^*$, $A_\zeta=A^*$, the following relations between the entries of the first moment $H_0(x)=H_0^*(x)$ hold
\begin{eqnarray*}
 \Re\pi_{12} =  -\dfrac{1}{2} \dfrac{\tau''}{\tau}, \quad
 \pi_{22} =  \dfrac{1}{3} \dfrac{\tau'''}{\tau}, \quad
 \Re \pi_{13}  =  -\dfrac{1}{6} \dfrac{\tau'''}{\tau}, \\
\Im \pi_{13}  = \pi_{11} \Im \pi_{12} + \Im \pi_{12}', \\
 \Re \pi_{23} = \dfrac{\Im \pi_{12}^2}{2} + \dfrac{2}{9} (q(x)^2-\dfrac{1}{4} q''(x)) + \dfrac{1}{8} \dfrac{\tau^{(4)}}{\tau}
\end{eqnarray*}
\end{lemma}

\section{\label{sec:GenVessels}Vessels with unbounded operators. Standard construction of a vessel}
The ideas presented in this Section can be found in \cite{bib:UnboundedVessels} for the symmetric case.
The class of functions serving as ``initial conditions'' for the transfer functions of vessels is defined as follows
\begin{defn} Class $\R(\sigma_1)$ consist of $p\times p$ matrix-valued functions $S(\lambda)$ of the complex variable $\lambda$, possessing the following representation:
\begin{equation} \label{eq:S0realized}
S(\lambda) = I - C_0 \mathbb X_0^{-1} (\lambda I - A)^{-1} B_0 \sigma_1
\end{equation}
where for an auxiliary Hilbert space $\mathcal H$ there are defined operators  $C_0:\mathcal H\rightarrow\mathbb C^3$,
$A_\zeta, \mathbb X_0, A: \mathcal H\rightarrow\mathcal H$,  $B_0:\mathbb C^3\rightarrow\mathcal H$. A general matrix-function $S(\lambda)$,
representable in such a form is called \textbf{realized}. Moreover, the operators are subject to the following assumptions:
\begin{enumerate}
	\item the operators $A, A_\zeta$ have dense domains $D(A), D(A_\zeta)$. $A, A_\zeta$ are generators of $C_0$ semi-groups on $\mathcal H$. 
	Denote the resolvents as follows $R(\lambda)=(\lambda I - A)^{-1}$, $R_\zeta(\lambda) = (\lambda I + A_\zeta)^{-1}$,
	\item the operator $B_0$ satisfies $R(\lambda) B_0 e \in \mathcal H$ for all $\lambda\not\in\SPEC(A), e\in\mathbb C^3$,
	\item the operator $\mathbb X_0$ is bounded and invertible,
	\item the \textbf{Lyapunov equation} holds for all $\lambda\not\in\SPEC(A)\cup\SPEC(-A_\zeta)$:
	\begin{equation} \label{eq:LyapunovAt0}
	 R(\lambda) \mathbb X_0 - \mathbb X_0 R_\zeta(\lambda) + R(\lambda) B_0 \sigma_1 C_0 R_\zeta(\lambda) = 0.
	\end{equation}
\end{enumerate}
We call an element of $\R(\sigma_1)$ as scattering matrix-function.
The subclass $\U(\sigma_1)\subseteq\R(\sigma_1)$ consists of symmetric functions, i.e. satisfying $S(\lambda)\sigma_1^{-1}S^*(-\lambda^*)=\sigma_1^{-1}$.
The Schur class $\SU(\sigma_1)\subseteq\U$ consists of symmetric functions, for which $\mathbb X_0$ is a positive operator. The sub-classes of rational functions
in $\SU, \U, \R$ are denoted by $\RSU, \RU, \RR$ respectively.
\end{defn}
When $S(\lambda)$ is just analytic at infinity (hence $A$ must be bounded), 
there is a very well known theory of realizations developed in \cite{bib:bgr}. 
For analytic at infinity and symmetric, i.e. satisfying $S^*(-\bar\lambda) \sigma_1 S(\lambda) =  \sigma_1$,
functions there exists a good realization theory using Krein spaces ($\mathcal H$ is a Krein space), 
developed in \cite{bib:KreinReal}\footnote{At the paper \cite{bib:KreinReal} a similar result is proved for functions symmetric with respect to the unit circle,
but it can be translated using Calley transform into $S^*(-\bar\lambda) \sigma_1 S(\lambda) =  \sigma_1$ and was done in \cite{bib:GenVessel, bib:SchurIEOT}}.
Such a realization is then translated into a function in $\U(\sigma_1)$.
The sub-classes $\U, \SU$ appear a lot in the literature and correspond to the symmetric case. We will not consider these two classes in this work and refer to \cite{bib:SchurIEOT}.

Now we present the \textit{standard construction} of a vessel $\mathfrak V$ from the
given realized matrix-function $S(\lambda)$ \eqref{eq:S0realized}:
\begin{enumerate}
	\item Let $B(x)$ be the unique solution of the following equation ($R(\lambda)=(\lambda I - A)^{-1}$)
	\begin{equation} \label{eq:RDB}
	\frac{\partial}{\partial x} R(\lambda) B  = - (A \,R(\lambda) B \sigma_2 + R(\lambda) B \gamma)\sigma_1^{-1}
	\end{equation}
	satisfying $B(x_0)=B_0$. 
	This equation is solvable because the coefficients $\sigma_1, \sigma_2, \gamma$ are constant and $A$ is
	a generator of a $C_0$ semigroup. For this equation to 
	hold we must require the following regularity assumptions:
	\begin{eqnarray}
	\label{eq:ResBsigma2} \forall\lambda\not\in\SPEC(A): R(\lambda) B(x)\sigma_2\mathbb C^3\subseteq D(A), \\ 
	\label{eq:ResBgamma} \forall\lambda\not\in\SPEC(A): R(\lambda) B(x)\gamma\mathbb C^3\subseteq\mathcal H,
	\end{eqnarray}
	\item Let $C(x)$ be the unique solution of \eqref{eq:DC} with the initial condition $C(x_0)=C_0$, 
	when we consider the equation \eqref{eq:DC}, applied to vectors in the dense set $D(A_\zeta)$ only,
	\item Solve \eqref{eq:DX} for $\mathbb X(x)$, satisfying $\mathbb X(x_0)=\mathbb X_0$ and let $\mathrm I$ be an interval, including $x_0$ on which $\mathbb X(x)$
	is invertible. Notice that $B(x)\sigma_2\mathbb C^3\subseteq\mathcal H$ follows from \eqref{eq:ResBsigma2}.
	\item Define $\gamma_*(x)$ on $\mathrm I$ by \eqref{eq:Linkage}.
\end{enumerate}
The main reason, why we call $S(\lambda)$ as the ``scattering data'' is the fact that $\gamma_*(x)$ (generalized
potential) is uniquely determined from $S(\lambda)$ by this construction. 
The question of uniqueness of $S(\lambda)$ for a given potential $\gamma_*(x)$ will be studied in Section \ref{sec:Uniqeness}.
A vessel in a more general form is defined as follows
\begin{defn} \label{def:VesselGen} The collection of operators and spaces 
\begin{equation} \label{eq:DefVGen}
\mathfrak{V} = (C(x), A_\zeta, \mathbb X(x), A, B(x); \sigma_1, \sigma_2, \gamma, \gamma_*(x);
\mathcal{H},\mathbb C^p; \mathrm I),
\end{equation}
is called a \textbf{vessel}, if $C(x):\mathbb C^3\rightarrow\mathcal H$, $A_\zeta, \mathbb X(x), A:\mathcal H\rightarrow\mathcal H$, 
$B(x):\mathbb C^p\rightarrow\mathcal H$ are differentiable linear operators, subject to the following conditions:
\begin{enumerate}
	\item the operators $A, A_\zeta$ have dense domains $D(A), D(A_\zeta)$. $A, A_\zeta$ are generators of $C_0$ semi-groups on $\mathcal H$. 
	Denote the resolvents as follows $R(\lambda)=(\lambda I - A)^{-1}$, $R_\zeta(\lambda) = (\lambda I + A_\zeta)^{-1}$,
	\item $B(x)$ satisfies regularity assumptions \eqref{eq:ResBsigma2}, \eqref{eq:ResBgamma}, and the equation \eqref{eq:RDB}.
	\item $C(x)$ satisfies \eqref{eq:DC} on $D(A_\zeta)$,
	\item $\mathbb X(x)$ is bounded, and invertible on $\mathrm I$ and satisfies \eqref{eq:DX},
	\item the \textbf{Lyapunov equation} holds for all $x\in\mathrm I, \lambda\not\in\SPEC(A)\cup\SPEC(-A_\zeta)$:
	\begin{equation} \label{eq:LyapunovGen}
		 R(\lambda) \mathbb X(x) - \mathbb X(x) R_\zeta(\lambda) + R(\lambda) B(x) \sigma_1 C(x) R_\zeta(\lambda) = 0.
	\end{equation}
	\item $\gamma_*(x)$ satisfies \eqref{eq:Linkage}.
\end{enumerate}
\end{defn}
The class of the transfer functions is defined as follows
\begin{defn} Class $\I=\I(\sigma_1,\sigma_1,\sigma_2,\gamma;I)$ 
consist of $3\times 3$ matrix-valued (transfer) functions $S(\lambda,x)$ of the complex variable $\lambda$ and $x\in\mathrm I$ for an interval $\mathrm I=[a,b]$ possessing the following representation:
\begin{equation} \label{eq:Srealized}
 S(\lambda,x) = I - C(x) \mathbb X^{-1}(x) (\lambda I - A)^{-1} B(x) \sigma_1 
\end{equation}
where for an auxiliary Hilbert space $\mathcal H$, the operators $C(x):\mathcal H\rightarrow\mathbb C^3$,
$A_\zeta, \mathbb X(x), A: \mathcal H\rightarrow\mathcal H$ and 
$B(x):\mathbb C^p\rightarrow\mathcal H$ constitute a vessel $\mathfrak V$ \eqref{eq:DefVGen} for some $\sigma_2,\gamma$.
\end{defn}
An analogues of the B\" acklund transformation Theorem \ref{thm:Backlund} for this new setting can be found in \cite{bib:UnboundedVessels} for the symmetric case. 
Actually, the assumptions on the operators were found in such a manner that  this theorem still holds. We omit its proof,
since we are interested in solutions of \eqref{eq:Boussinesq}.

\section{\label{sec:Uniqeness}``Uniqueness'' of the scattering data}
Let us consider now the uniqueness of the scattering matrix $S(\lambda,0)$. First we prove the following
\begin{lemma} Suppose that we are given a regular vessel \eqref{eq:DefV}
\[ \mathfrak{V}_{reg} = (C(x), A, \mathbb X(x), A_\zeta, B(x); \sigma_1, 
\sigma_2, \gamma, \gamma_*(x);
\mathcal{H},\mathbb C^3;\mathrm I), \]
realizing coefficients $q(x), p(x)$. Let $S(\lambda,x)$ be its transfer function, defined in \eqref{eq:DefS}.
Let $Y(\lambda)$ be an arbitrary $3\times 3$ matrix function, commuting with the fundamental solution 
$\Phi(\lambda,x)$ of \eqref{eq:InCC}. Then $\widetilde S(\lambda,x) = S(\lambda,x)Y(\lambda)$ is the transfer function of a vessel realizing the same
coefficients $q(x), p(x)$.
\end{lemma}
\noindent\textbf{Proof:} By the definition it follows that
\[ S(\lambda,x) = \Phi_*(\lambda,x) S(\lambda,0) \Phi^{-1}(\lambda,x).
\]
So,
\[ \widetilde S(\lambda,x) = S(\lambda,x) Y(\lambda) = \Phi_*(\lambda,x) S(\lambda,0) \Phi^{-1}(\lambda,x)Y(\lambda) = \Phi_*(\lambda,x) S(\lambda,0) Y(\lambda) \Phi^{-1}(\lambda,x)
\]
and realizes the same coefficients $q(x), p(x)$. By the standard construction, there is a vessel $\widetilde{\mathfrak V}$, whose transfer functions is $\widetilde S(\lambda,x)$.
\qed

Let us investigate the structure of a matrix $Y(\lambda)$, commuting with $\Phi(\lambda,x)$.
Using the form \eqref{eq:PhiForm}, it is easy to conclude that a matrix, which commutes with $\Phi(\lambda,x)$ must be of the form
\begin{equation} \label{eq:YForm}
 Y(\lambda) = a(\lambda) I + b(\lambda) \bbmatrix{0&1&0\\0&0&-1\\-\lambda&0&0}+c(\lambda)\bbmatrix{0&0&1\\-\lambda&0&0\\0&\lambda&0},
\end{equation}
by considering the coefficients of pure exponents in $\Phi(\lambda,x)Y(\lambda) = Y(\lambda)\Phi(\lambda,x)$.
\begin{thm}[Uniquness of the scattering matrix]\label{thm:Uniqueness} Suppose that $S(\lambda,x)$, $\widetilde S(\lambda,x)$ are the transfer functions of two regular vessels $\mathfrak V_{reg}, \widetilde {\mathfrak V}_{reg}$,
defined in \eqref{eq:DefV}, realizing the same potential $\gamma_*(x)$. Then there exists a matrix $Y(\lambda)\in\R$ such that
\[ \widetilde S(\lambda,x) = S(\lambda,x) Y(\lambda).
\]
\end{thm}
\Proof Let us consider the function $S^{-1}(\lambda,x)\widetilde S(\lambda,x)$. By the definition this functions maps solutions of the input LDE \eqref{eq:InCC} to itself:
\begin{multline*}
 S^{-1}(\lambda,x)\widetilde S(\lambda,x) = \left(\Phi_*(\lambda,x) S(\lambda,0) \Phi^{-1}(\lambda,x)\right)^{-1}\Phi_*(\lambda,x) \widetilde S(\lambda,0) \Phi^{-1}(\lambda,x) =\\
= \Phi(\lambda,x)  S^{-1}(\lambda,0)\widetilde S(\lambda,0) \Phi^{-1}(\lambda,x)
\end{multline*}
Plug here, the formula \eqref{eq:PhiForm} and find conditions so that the coefficients of the exponents $e^{-kx}, e^{-\alpha k x}, e^{-\alpha^2 k x}$ vanish. This is necessary for
making this function bounded at infinity out of the spectrum of $A$. Then calculations show that actually $S^{-1}(\lambda,0)\widetilde S(\lambda,0)$ must commute with
$\Phi(\lambda,x)$ so that this functions cancels all the singularities at infinity. As a result, by the preceding arguments it must be a function $Y(\lambda)$ of the form
\eqref{eq:YForm}. And we obtain that
\[ S^{-1}(\lambda,0)\widetilde S(\lambda,0) = Y(\lambda),
\]
from where the result follows.
\qed

Another, weaker form of the uniqueness is used later in the text and is presented in the next Lemma. We emphasize that a similar theorem lemma was proved in 
the Sturm-Liouville case in \cite{bib:GenVessel} and in \cite{bib:FaddeyevII} for purely continuous spectrum.
\begin{lemma}\label{lemma:uniquegamma*} Suppose that two functions $S(\lambda,x)$, $\widetilde S(\lambda,x)$ are in class $\I(\sigma_1, \sigma_2,\gamma)$, possessing the same initial value
\[ S(\lambda,0) = \widetilde S(\lambda,0)
\]
and are bounded at a neighborhood of infinity, with a limit value $I$ there.
Then the corresponding outer potentials are equal:
\[ \gamma_*(x) = \widetilde \gamma_*(x).
\]
\end{lemma}
\Proof Suppose that
\[ S(\lambda,x) = \Phi_*(\lambda,x) S(\lambda,0) \Phi^{-1}(\lambda,x), \quad \widetilde S(\lambda,x) = \widetilde \Phi_*(\lambda,x) S(\lambda,0) \Phi^{-1}(\lambda,x).
\]
Then
\[ \widetilde S^{-1}(\lambda,x) S(\lambda,x) = \widetilde \Phi_*(\lambda,x)\Phi^{-1}_*(\lambda,x)
\]
is entire (the singularities appear in $S(\lambda,0) = \widetilde S(\lambda,0)$ only and are cancelled) and equal to $I$ (- the identity matrix) at infinity.
By a Liouville theorem, it is a constant function, namely $I$. So $\widetilde \Phi_*(\lambda,x)\Phi^{-1}_*(\lambda,x) = I$ or
\[ \widetilde \Phi_*(\lambda,x) = \Phi_*(\lambda,x).
\]
If we differentiate this, we obtain that $\widetilde \gamma_*(x)=\gamma_*(x)$.
\qed
\section{Choice of the parameters realizing the Boussinesq equation \eqref{eq:Boussinesq}}
\subsection{Realizing the Boussinesq equation by a Boussinesq vessel}
Let us choose the following parameters
\begin{equation} \label{eq:3DimCanParamt}
\widetilde{\sigma}_1 = \sigma_1, \quad \widetilde\sigma_2=\bbmatrix{0&-i&0\\i&0&0\\0&0&0},\quad \widetilde\gamma=\bbmatrix{0&0&0\\0&0&0\\0&0&i}.
\end{equation}
Suppose also that the operators $C,\mathbb X,B$ evolve with respect to $t$ using the same formulas as for the vessel with vessel parameters 
$\sigma_1,\sigma_2,\gamma$ substituted with $\widetilde{\sigma}_1, \widetilde\sigma_2, \widetilde\gamma$. More precisely, we define the following vessel:
\begin{defn} Suppose that the parameters $\sigma_1,\sigma_2,\gamma$ are defined in Definition \ref{def:3DimCanParam} and
$\widetilde{\sigma}_1, \widetilde\sigma_2, \widetilde\gamma$ are defined in \eqref{eq:3DimCanParamt}. Then a (regular non-symmetric) Boussinesq vessel $\mathfrak V_{Bouss,reg}$ is a collection
of operators and spaces and a rectangle $\mathrm R$
\begin{equation} \label{eq:DefVBouss}
\mathfrak{V}_{Bouss, reg} = (C(x,t), A_\zeta, \mathbb X(x,t), A, B(x,t); \sigma_1, 
\sigma_2, \gamma, \gamma_*(x,t),\widetilde{\sigma}_1, \widetilde\sigma_2, \widetilde\gamma ;
\mathcal{H},\mathbb C^3;\mathrm R),
\end{equation}
where the bounded operators $C(x,t):\mathcal H\rightarrow\mathbb C^3$, $A_\zeta,\mathbb X(x,t),A:\mathcal H\rightarrow\mathcal H$, $B(x,t):\mathbb C^3\rightarrow\mathcal H$ 
and a $3\times 3$ matrix function $\gamma_*(x,t)$ satisfy the vessel conditions \eqref{eq:DB}, \eqref{eq:DC}, \eqref{eq:DX}, \eqref{eq:Linkage}, \eqref{eq:Lyapunov}
and the following evolutionary equations
\begin{eqnarray}
\label{eq:DBt} \frac{\partial}{\partial t} B  &  & = - (A \, B \widetilde\sigma_2 + B \widetilde\gamma)\widetilde\sigma_1^{-1}, \\
\label{eq:DCt} \frac{\partial}{\partial t} C  &  & = \widetilde\sigma_1^{-1} (\widetilde\gamma C - \widetilde\sigma_2 C A_\zeta ), \\
\label{eq:DXt} \frac{\partial}{\partial t} \mathbb X &  & =  B \widetilde\sigma_2 C.
\end{eqnarray}
the operator $\mathbb X(x,t)$ is assumed to be invertible on the rectangle $\mathrm R$. In the case $A_\zeta=A^*$ and $C=B^*$, we call such a vessel \textbf{symmetric}.
\end{defn}
\begin{thm}\label{thm:BoussnesqProof} Suppose that $\mathfrak{V}_{Bouss,reg}$ is a Boussinesq regular vessel, defined in \eqref{eq:DefVBouss} and $\tau(x)=\det(\mathbb X(x_0,t_0)^{-1}\mathbb X(x,t))$ is its tau function,
defined for arbitrary point $(x_0,t_0)\in\mathrm R$.
Then the coefficient $q(x)=-\dfrac{3}{2} \dfrac{\partial^2}{\partial x^2} \ln (\tau(x,t))$ satisfies the Boussinesq equation \eqref{eq:Boussinesq} on $\mathrm R$.
\end{thm}
\noindent\textbf{Proof:} Consider the Boussinesq equation \eqref{eq:Boussinesq}:
\[
q_{tt} = \dfrac{\partial^2}{\partial x^2} [q_{xx} - 4 q^2 ].
\]
From the formula \eqref{eq:qtau} we find that it can be rewritten as
\[
-\dfrac{3}{2} (\pi_{11})_{xtt} = \dfrac{\partial^2}{\partial x^2} [-\dfrac{3}{2} \pi_{xxx}- 4 (-\dfrac{3}{2} (\pi_{11})_x)^2 ] 
\]
or integrating with respect to $x$ once and multiplying by $-\dfrac{2}{3}$ it is enough to show that:
\begin{equation} \label{eq:enoghtoprove}
(\pi_{11})_{tt} = (\pi_{11})_{xxxx} + 12 (\pi_{11})_x (\pi_{11})_{xx}.
\end{equation}
The formula \eqref{eq:DHn}, where we substitute the vessel parameters with tilde notation, and substitute 
$x$-derivative with the $t$-derivative, implies that
\begin{equation} \label{eq:DHnt}
 \widetilde\sigma_1^{-1}\widetilde\sigma_2H_{n+1} - H_{n+1} \widetilde\sigma_2\widetilde\sigma_1^{-1} = (H_n)'_t - 
\widetilde\sigma_1^{-1} (\widetilde\gamma+\widetilde\sigma_2H_0\widetilde\sigma_1-\widetilde\sigma_1H_0\widetilde\sigma_2)  H_n + H_n \widetilde\gamma\widetilde\sigma_1^{-1}.
\end{equation}
Since $\pi_{11}=\TR(\sigma_2 H_0)$, we can use \eqref{eq:DHnt} to obtain:
\[ \dfrac{\partial}{\partial t} \pi_{11} = \dfrac{\partial}{\partial t} \TR(\sigma_2 H_0) =
\TR(\sigma_2 \dfrac{\partial}{\partial t} H_0) = i (\pi_{12}'-\pi_{21}') = 
\dfrac{\partial}{\partial x} \TR(\bbmatrix{0&-i&0\\i&0&0\\0&0&0}H_0),
\]
after cancellations and using \eqref{eq:pi21}. So, differentiating again and using \eqref{eq:DHnt} we obtain:
\[ \begin{array}{llll}
\dfrac{\partial^2}{\partial t^2} \pi_{11} & = \dfrac{\partial}{\partial t} \dfrac{\partial}{\partial x} \TR(\bbmatrix{0&-i&0\\i&0&0\\0&0&0}H_0) =
\dfrac{\partial}{\partial x} \TR(\bbmatrix{0&-i&0\\i&0&0\\0&0&0}\dfrac{\partial}{\partial t}H_0) = \\
& = - [\pi_{12}\pi_{11}']' -  \pi_{11}^2 \pi_{11}'' -\dfrac{1}{2} \pi_{11}'\pi_{11}''+\dfrac{1}{2}\pi_{11}(-4(\pi_{11}')^2 + \pi_{11}''') + \dfrac{1}{2}\pi_{11}''''.
\end{array} \]
Plugging here \eqref{eq:pi12relation} we obtain after cancellations \eqref{eq:enoghtoprove}. \qed

\subsection{Standard construction of a Boussinesq vessel}
We may revise now the standard construction of a vessel, presented in Section \ref{sec:GenVessels} in order to show how to construct a Boussinesq vessel
and hence how to produce a solution of the Boussinesq equation \eqref{eq:Boussinesq}. For this we have to continue the process of construction, using the parameters
$\widetilde{\sigma}_1, \widetilde\sigma_2, \widetilde\gamma$ defined in \eqref{eq:3DimCanParamt}:
\begin{enumerate}
	\item[1-3] use old steps to construct $C(x)$, $\mathbb X(x)$, $B(x)$ with the same assumptions, as in Section \ref{sec:GenVessels},
	\item[4.] solve the equation
	\begin{equation} \label{eq:RDBt}
	\frac{\partial}{\partial x} R(\lambda) B  = - (A \,R(\lambda) B \widetilde\sigma_2 + R(\lambda) B \widetilde\gamma)\widetilde\sigma_1^{-1}
	\end{equation}
	with the initial condition $B(x)$ resulting in a function $B(x,t)$.
	We have to require, similarly to the construction of $B(x)$ that
	\begin{eqnarray}
	\label{eq:ResBsigma2t} \forall\lambda\not\in\SPEC(A): R(\lambda) B(x,t)\widetilde\sigma_2\mathbb C^3\subseteq D(A), \\ 
	\label{eq:ResBgammat} \forall\lambda\not\in\SPEC(A): R(\lambda) B(x,t)\widetilde\gamma\mathbb C^3\subseteq\mathcal H,
	\end{eqnarray}
	\item solve for $C(x,t)$ the equation \eqref{eq:DCt} with initial condition $C(x)$ on $D(A_\zeta)$,
	\item solve for $\mathbb X(x,t)$ the equation \eqref{eq:DXt} with initial $\mathbb X(x)$,
	\item define $\gamma_*(x,t)$ by \eqref{eq:Linkage} at all points, where $\mathbb X(x,t)$ is invertible.
\end{enumerate}
Meanwhile, we present a weaker, regular case form of a general Theorem \ref{thm:BoussSolGen} on producing
solutions of \eqref{eq:Boussinesq}:
\begin{thm} Suppose that a collection
\[ (C(x,t), A_\zeta, \mathbb X(x,t), A, B(x,t); \sigma_1, 
\sigma_2, \gamma, \gamma_*(x,t),\widetilde{\sigma}_1, \widetilde\sigma_2, \widetilde\gamma ;
\mathcal{H},\mathbb C^3;\mathrm R)
\]
is obtained from a regular (i.e. all the initial operators are bounded) $S(\lambda)\in\R$ by the standard construction around a point $(x_0,t_0)$. 
Then this collection is a Boussinesq vessel on a rectangle $\mathrm R$, including $(x_0,t_0)$. The coefficient
$q(x) = -\dfrac{3}{2} \dfrac{\partial^2}{\partial x^2} (\mathbb X^{-1}(x_0,t_0)\mathbb X(x,t))$ satisfies the Boussinesq equation \eqref{eq:Boussinesq} on $\mathrm R$.
\end{thm}
\Proof Notice that since $\mathbb X(x_0,t_2)$ is invertible, by a continuity, since all the operators are bounded, it will be invertible on a rectangle 
$\mathrm R$, including $(x_0,t_0)$. Then the Theorem follows from the definitions and Theorem \ref{thm:BoussnesqProof}. \qed
\section{Examples of solutions of the Boussinesq equation \eqref{eq:Boussinesq}}
\subsection{\label{sec:solitons}Solitons}
If we take $\mathcal H = \mathbb C$, we will obtain solitons. Use the following choice of
solutions of the vessel equations
\[ \begin{array}{llll}
A =  (2 i \mu)^3, \quad  A\zeta = (i \mu)^3, \\
B(x,t) = \bbmatrix{E_B &-2i \alpha\mu E_B &(2\alpha\mu)^2 E_B },  & E_B = \exp(2i\alpha\mu(x+2\alpha\mu t))\\
C(x,t) = \bbmatrix{E_C \\-i \alpha\mu E_C\\(\alpha\mu)^2E_C}, & E_C= \exp(i\alpha\mu(x-\alpha\mu t)), \\
\mathbb X(x,t) = -\dfrac{2\alpha\cosh(\dfrac{\sqrt{3}}{2} \mu(x + \mu t)) E_X }{\sqrt{3}\mu},
& E_X = \exp\left(\dfrac{\sqrt{3}}{2} \mu x(2\alpha^2-1) + \dfrac{\sqrt{3}}{2} \mu^2 t (2\alpha^2+5)\right),
\end{array} \]
where $\alpha=e^{\frac{2i\pi}{3}}$ is the basic third root of $1$.
Then simple calculations show that
\begin{equation} \label{eq:classSol}
 q(x) = -\dfrac{9\mu^2}{2\cosh^2(\dfrac{\sqrt{3}}{2} \mu(x+t\mu))}
\end{equation}
is the Boussinesq \textit{soliton} solution of \eqref{eq:Boussinesq} (see \cite{bib:SolitonBouss}).

Another non-trivial soliton, obtained in this form is as follows.
\[ \begin{array}{llll}
A =  (2 i \mu)^3, \quad  A\zeta = (2i \mu)^3, \quad k_1 = \alpha 2 i \mu, \quad f_1 = \alpha^2 2 i \mu, \\
B(x,t) = \bbmatrix{b \dfrac{E_B}{k_1^2} + c \dfrac{E_f}{f_1^2} & - b \dfrac{E_B}{k_1} - c \dfrac{E_f}{f_1} &-b E_B - c E_f},  & E_B = \exp(k_1 x- i k_1^2 t)\\
& E_f = \exp(f_1 x- i f_1^2 t), \\  c=-i, \quad b= \dfrac{1}{i+\sqrt{3}} \\
C(x,t) = \bbmatrix{\dfrac{E_f}{f_1^2} \\-\dfrac{E_f}{f_1}\\-E_f},  \\
\mathbb X(x,t) = \dfrac{\alpha e^{-2ix\mu}(e^{2\sqrt{3}x\mu} + e^{4\sqrt{3}t\mu^2})}{64 \mu^5}
\end{array} \]
and the coefficient $q(x,t)$ is
\begin{equation} \label{eq:solnonclass}
 q(x,t) = -\dfrac{18e^{2\sqrt{3}\mu(x+2t\mu)}\mu^2}{(e^{2\sqrt{3}x\mu} + e^{4\sqrt{3}t\mu^2})^2}
\end{equation}
\subsection{Solutions, belonging to the Schwartz class.}
Suppose that we are given real coefficients $q(x), p(x)$ in the Schwartz class and define $\pi_{11}$ by formula \eqref{eq:qtau}.
Suppose that the functions $F_1, F_2, F_3$ are solutions of \eqref{eq:Lisk3} satisfying
\[ \begin{array}{llll}
F_1(x,k) = e^{-kx}(1 + O(1)), &\text{ as } &x\rightarrow+\infty, \\
F_2(x,k) = e^{-\alpha k x}(1 + O(1)),& \text{ as } &x\rightarrow+\infty, \\
F_3(x,k) = e^{-\alpha^2 k x}(1 + O(1)),& \text{ as } &x\rightarrow+\infty 
\end{array} \]
where $k\in\Omega_3=\{ k| -\dfrac{2\pi}{3}<\arg(k)<\dfrac{2\pi}{3}\}$.
For the coefficients $q(x), p(x)\in\mathcal S$, where $\mathcal S$ - the Schwartz class of rapidly decreasing functions, the existence of 
$F_1, F_2, F_3$ is shown in \cite{bib:DTTBoussinesq}, for example.

Then a function-matrix $\Psi_*(\lambda,x)$ which is a solution of \eqref{eq:OutCC} can be constructed from 
solutions $\Psi_{*,1}, \Psi_{*,2}, \Psi_{*,3}$ of \eqref{eq:Lisk3} and is given in view of \eqref{eq:entries} by
($k=\sqrt[3]{\lambda}$)
\begin{multline} \Psi_*(\lambda,x) = \\
\dfrac{1}{3\alpha^2} \bbmatrix{
\Psi_{*,1} & \Psi_{*,2} & \Psi_{*,3} \\
-\pi_{11}\Psi_{*,1}-\Psi_{*,1}' & -\pi_{11}\Psi_{*,2}-\Psi_{*,2}' & -\pi_{11}\Psi_{*,3}-\Psi_{*,3}' \\
\pi_{21}\Psi_{*,1} -(\pi_{11}\Psi_{*,1}+\Psi_{*,1}')'  & \pi_{21}\Psi_{*,2} -(\pi_{11}\Psi_{*,2}+\Psi_{*,2}')'  & \pi_{21}\Psi_{*,3} -(\pi_{11}\Psi_{*,3}+\Psi_{*,3}')' 
}
\end{multline}
where 
\[ \Psi_{*,1} = \alpha^2(F_1 + F_2 + F_3), \quad
\Psi_{*,2} = \dfrac{\alpha^2}{k}F_1 +\dfrac{\alpha}{k} F_2+\dfrac{1}{k}F_3, \quad
\Psi_{*,3} = \dfrac{\alpha^2}{k^2} F_1+\dfrac{1}{k^2} F_2+\dfrac{\alpha}{k^2} F_3.
\]
As a function of $\lambda$ the matrix-function $\Psi_*(\lambda,x)$ is analytic in the whole plane except for the cut along the negative axis $(-\infty,0]$.
Notice that this function is globally bounded as there are asymptotic formulas for such solutions appearing in \cite{bib:DTTBoussinesq}. Actually, the three functions 
$F_1, F_2, F_3$ have the asymptotic behavior $\exp(\alpha^i kx) (1 + O(k)), i=0,1,2$ as $k\rightarrow\infty$ for all $x\in\mathbb R$ and these asymptotics can be differentiated infinitely
many times.

So we define 
\[ S(\lambda,x) = \Psi_*(\lambda,x) \Phi^{-1}(\lambda,x),
\]
which becomes
\[ S(\lambda,x) = \Phi_*(x,\lambda) S(\lambda,0)\Phi^{-1}(\lambda,x)
\]
with $\Phi_*(x,\lambda)$ - the fundamental solution of \eqref{eq:OutCC} equal to identity at $x=0$. Since
$\Psi_*$ and $\Phi$ are analytic functions of $k$ in $\Omega_3$, they will be analytic functions of $\lambda$,
except for the cut along the negative real axis, where it will have jumps, but is remained bounded at the infinity.
These jumps are expressible by parameters $R_1, R_2$, appearing at \cite{bib:DTTBoussinesq} in a quite
general case \cite[Theorem 14]{bib:DTTBoussinesq} of real coefficients $q(x), p(x)$.

Applying the standard construction to 
$S(\lambda,0)$ we will obtain a vessel realizing $q_V(x)$ and derived from it $p_V(x)$, which is equal by Lemma \ref{lemma:uniquegamma*} to the given $q(x)$: $q_V(x)=q(x)$. Evolving further the operators with respect to $t$,
we will create a solution of \eqref{eq:Boussinesq}, coinciding with $q(x)$ at time $t=0$.

\subsection{General solutions}
We can obtain more general solutions of \eqref{eq:Boussinesq} if we choose to work in a more general setting, presented in Section \ref{sec:GenVessels}.
Applying the standard vessel construction to a function $S(\lambda)\in\R$, we obtain 
\begin{thm} \label{thm:BoussSolGen}
Suppose that a Boussinesq vessel $\mathfrak V_{Bouss}$ is obtained from 
\[ S(\lambda) = I - C_0\mathbb X^{-1}(\lambda I - A)^{-1}B_0\sigma_1
\]
by the standard construction, satisfies the following assumptions ($\forall (x,t)\in\mathbb R^2$):
\begin{eqnarray}
\label{eq:BregBouss} B(x,t) \bbmatrix{1\\0\\0}\in D(A), \quad B(x,t)\bbmatrix{0\\1\\0}\in D(A) , \quad B(x,t)\bbmatrix{0\\0\\1} \in\mathcal H.
\end{eqnarray}
Then the coefficient $q(x,t)$ is twice differentiable with respect to $t$ and four times differentiable with respect to $x$ on an open subset 
$\Omega \subseteq (\mathbb R\times\mathbb R)\backslash Z$, obtained after removing a closed subset $Z\subseteq\mathbb R\times\mathbb R$ of the points 
$(x,t)\in Z$, in which $\mathbb X(x,t)$ is not invertible. Moreover, the coefficient $q(x,t)$ satisfies the Boussinesq equation \eqref{eq:Boussinesq} on $\Omega$.
\end{thm}
\Rems Notice that in the regular case the operator $\mathbb X(x,t)$ and hence the tau function 
$\tau(x,t)=\det(\mathbb X^{-1}(x_0,t_0)\mathbb X(x,t))$ are analytic
on $\mathbb R\times\mathbb R$. So, the zeros of $\tau(x,t)$ are analytic sets of $\mathbb R\times\mathbb R$, but in general their structure is complicated and out of the scope of this work. \textbf{2.} The assumptions \eqref{eq:BregBouss} imply the regularity assumptions \eqref{eq:ResBsigma2},	\eqref{eq:ResBgamma}, \eqref{eq:ResBsigma2t},	\eqref{eq:ResBgammat}.

\Proof Notice that the regularity assumptions \eqref{eq:ResBsigma2},	\eqref{eq:ResBgamma} and \eqref{eq:ResBsigma2t},	\eqref{eq:ResBgammat}
imply that the operators $B(x,t):\mathbb C^3\rightarrow\mathcal H$, $C(x,t):\mathbb C^3\rightarrow D(A_\zeta)$ are defined for all $(x,t)\in\mathbb R^2$ (actually, for 
$C(x) = \bbmatrix{c_1(x,t)&c_2(x,t)&c_3(x,t)}$ one has to verify that the three entries $c_1(x,t), c_2(x,t), c_3(x,t)$ are in $D(A_\zeta)$ for all $x$).
The operator $\mathbb X(x,t)$ is the unique solution of \eqref{eq:DX}, \eqref{eq:DXt} with the initial condition $\mathbb X_0$. Notice that this is an operator of the form
$\mathbb X_0 + T(x,t)$ for a trace-class operator $T(x,t)$ and we can define $\tau(x,t)=\det(\mathbb X_0^{-1}\mathbb X(x,t))$ \eqref{eq:Deftau}.
Moreover, since $\mathbb X(x)$ is $x$-differentiable
\[  \dfrac{\tau'(x,t)}{\tau(x,t)} = \TR(\mathbb X'(x,t)\mathbb X^{-1}(x,t)) = \TR(B(x,t)\sigma_2C(x)\mathbb X^{-1}(x,t))
\]
for all points $(x,t)\in Z=\{ (x,t)\mid \tau(x,t)\neq 0 \}) = \{ (x,t)\mid \mathbb X(x,t) \text{ is not invertible}\}$. On the other hand, the
differentiability of the coefficient $q(x) = -\dfrac{3}{2} \dfrac{\partial^2}{\partial x^2} \ln \tau(x,t)$ follows from the differentiability of $\dfrac{\tau'(x,t)}{\tau(x,t)}$.
So, we investigate the existence of the derivatives for $\dfrac{\tau'(x,t)}{\tau(x,t)}$.

Applying the formulas \eqref{eq:DB}, \eqref{eq:DCX} we obtain that
\begin{multline*}
 \dfrac{\tau''(x,t)}{\tau(x,t)} =  \dfrac{\partial}{\partial x} \dfrac{\tau'(x,t)}{\tau(x,t)} + (\dfrac{\tau'(x,t)}{\tau(x,t)})^2 = \TR(B(x,t)(\sigma_2\sigma_1^{-1}\gamma-\gamma\sigma_1^{-1}\sigma_2) C(x,t)\mathbb X^{-1}(x,t)) = \\
 = - \TR(\bbmatrix{0&1&0\\1&0&0\\0&0&0} C(x,t)\mathbb X^{-1}(x,t)B(x)),
\end{multline*}
which is well-defined on $\Omega$. Then, similarly,
\[
\dfrac{\tau'''(x,t)}{\tau(x,t)} =  \dfrac{\partial}{\partial x} \dfrac{\tau''(x,t)}{\tau(x,t)} + \dfrac{\tau''(x,t)}{\tau(x,t)} \dfrac{\tau'(x,t)}{\tau(x,t)}
= -\TR(\bbmatrix{0&0&1\\0&-2&0\\1&0&0} C(x,t)\mathbb X^{-1}(x,t)B(x)).
\]
Then, differentiating and plugging \eqref{eq:DB}, \eqref{eq:DCX} again we can obtain that
\[
\dfrac{\partial}{\partial x} \dfrac{\tau'''(x,t)}{\tau(x,t)} = 3 (\pi_{12}\pi_{21} - \pi_{23}-\pi_{32}).
\]
Since $\pi_{12}=\bbmatrix{1&0&0}C(x,t)\mathbb X^{-1}(x,t)B(x,t)\bbmatrix{0\\1\\0}$ in view of \eqref{eq:ResBgamma} it is a well-defined function on $\Omega$.
Similarly, $\pi_{21}, \pi_{23}, \pi_{32}$ are well defined from \eqref{eq:ResBgamma}.
Finally, differentiating the last expression, we will obtain that
\[ \dfrac{\partial^2}{\partial x^2} \dfrac{\tau'''(x,t)}{\tau(x,t)} = 3 (\pi'_{12}\pi_{21} + \pi_{12}\pi'_{21} - \pi'_{23}-\pi'_{32})
\]
Here, the derivative of $\pi_{12}$ is as follows
\begin{multline*}
 \pi'_{12} = \dfrac{\partial}{\partial x} \left(\bbmatrix{1&0&0}C(x,t)\mathbb X^{-1}(x,t)B(x,t)\bbmatrix{0\\1\\0}\right) = \\
= \bbmatrix{1&0&0} \sigma_1 \gamma_*(x,t) C(x,t)\mathbb X^{-1}(x,t)B(x,t)\bbmatrix{0\\1\\0} -
\bbmatrix{1&0&0} C(x,t)\mathbb X^{-1}(x,t)B(x,t) \gamma \sigma_1 \bbmatrix{0\\1\\0},
\end{multline*}
where we used \eqref{eq:DB}, \eqref{eq:DCX}. Thus it is a well-defined function on $\Omega$. Similarly for $\pi'_{21}$. The derivatives of $\pi_{23}, \pi_{32}$ are as follows
\[ \begin{array}{lll}
\pi'_{32} = - \bbmatrix{0&1&0}C(x,t)\mathbb X^{-1}(x,t) A B(x,t)\bbmatrix{1\\0\\0} - \pi_{12} \pi_{13} + \pi_{33}, \\
\pi'_{23} = \bbmatrix{1&0&0}C(x,t)\mathbb X^{-1}(x,t) A B(x,t)\bbmatrix{0\\1\\0} + \pi_{12}(\pi_{13} + \pi_{22}
 - \pi_{31}) + \pi_{11} \pi_{32} + \pi_{33},
\end{array} \]
which are well defined by the assumption \eqref{eq:BregBouss} on $\Omega$.

The fact that we can $t$-differentiate $q(x,t)$ twice follows exactly the same lines, using the parameters $\widetilde \sigma_1, \widetilde \sigma_2, \widetilde \gamma$ and the regularity assumption \eqref{eq:BregBouss}. Using \eqref{eq:DHnt}, the first derivative is
\[ \dot\pi_{11} = i (\pi_{13}-\pi_{31} + \pi_{11}(\pi_{21}-\pi_{12})),
\]
which is well-defined on $\Omega$. For the second derivatives, it is enough to check that $\dot\pi_{13}-\dot\pi_{31}$, $\dot\pi_{21}-\dot\pi_{12}$ are
well defined. Using \eqref{eq:DHnt} again
\[ \begin{array}{lll}
i \dot\pi_{13} =  \bbmatrix{0&1&0}C(x,t)\mathbb X^{-1}(x,t) A B(x,t)\bbmatrix{1\\0\\0} + \pi_{21} \pi_{22} + \pi_{11} (\pi_{23} + \pi_{32}) - \pi_{33}\\
-i \dot\pi_{12} = \bbmatrix{1&0&0}C(x,t)\mathbb X^{-1}(x,t) A B(x,t)\bbmatrix{1\\0\\0}  + \pi_{12} \pi_{21} + \pi_{21}^2 + \pi_{23} + 
 \pi_{11} (\pi_{13} - \pi_{22} + \pi_{31}),
\end{array} \]
which are well-defined on $\Omega$ by \eqref{eq:BregBouss}. Similarly for $\dot\pi_{21},\dot\pi_{31}$ and the proof is finished.
\qed

\Rem Notice that one more derivative of $\pi'_{23}$ produces a term of the form
\[ \bbmatrix{1&0&0}C(x,t)\mathbb X^{-1}(x,t) A B(x,t)\bbmatrix{0\\0\\1},
\]
which may fail to be a well-defined function, because $B(x,t)\bbmatrix{0\\0\\1}$ may fail to be at the domain of $A$.

\bibliographystyle{alpha}
\bibliography{../biblio}

\end{document}